\documentclass[11pt,reqno]{amsart}
\usepackage[letterpaper,textwidth=6.05in,textheight=8.60in,centering]{geometry}
\usepackage[T1]{fontenc}
\usepackage{lmodern,amssymb,mathtools,microtype}
\usepackage[hidelinks]{hyperref}
\newcommand{\R}{\mathbb R}
\newcommand{\T}{\mathbb T}
\newcommand{\dd}{\,\mathrm d}
\newcommand{\abs}[1]{\lvert#1\rvert}
\newcommand{\norm}[1]{\lVert#1\rVert}
\newtheorem{theorem}{Theorem}[section]
\newtheorem{lemma}[theorem]{Lemma}
\newtheorem{proposition}[theorem]{Proposition}
\newtheorem{corollary}[theorem]{Corollary}
\theoremstyle{remark}

\title{A Proof of the Ovals Conjecture}
\author{Durvudkhan Suragan}
\address{School of Computing and Artificial Intelligence, Nazarbayev University, Astana, Kazakhstan}
\email{durvudkhan.suragan@nu.edu.kz}
\date{}
\begin{document}
\begin{abstract}
We prove that the lowest eigenvalue of the periodic Schr\"odinger operator $-\partial_s^2+\kappa^2$ on a closed curve of length $\ell$ is at least $4\pi^2/\ell^2$. Equality holds precisely for the family of ovals described by Benguria and Loss. The proof uses a new comparison argument for scalar ordinary differential equations.
\end{abstract}
\maketitle

\section{Introduction and the main result}

Let $\gamma\colon\T_\ell\to\R^n$, $n\geq2$, be a closed curve of length $\ell$, parametrized by arclength, where $\T_\ell=\R/\ell\mathbb Z$. We assume $\gamma\in W^{2,2}(\T_\ell;\R^n)$ and write $T=\gamma'$ and $\kappa=\abs{T'}$. The first (lowest) eigenvalue of the periodic Schr\"odinger operator with potential $\kappa^2$ is
\begin{equation}\label{eq:rq}
 \lambda_1(\gamma)=\inf_{0\ne f\in H^1(\T_\ell)}
 \frac{\int_{\T_\ell} (\abs{f'}^2+\kappa^2\abs{f}^2)\dd s}
 {\int_{\T_\ell}\abs{f}^2\dd s}.
\end{equation}
The form is well defined since $\kappa^2\in L^1$ and $H^1(\T_\ell)$ embeds into $C(\T_\ell)$. Neither planarity nor simplicity is assumed.

Benguria and Loss \cite{BL} conjectured that the infimum of \eqref{eq:rq} over closed curves of length $2\pi$ is one. The circle attains this value, but it is by no means the only candidate for a minimizer. There is a family of noncircular ovals with the same eigenvalue which degenerates to a line segment traversed twice. This loss of compactness is one of the difficulties of the problem.

\begin{theorem}\label{thm:main}
For every closed curve as above,
\begin{equation}\label{eq:main}
 \lambda_1(\gamma)\geq \frac{4\pi^2}{\ell^2}.
\end{equation}
Writing $\omega=2\pi/\ell$, equality holds if and only if there are linearly independent vectors $A,B\in\R^n$ such that
\begin{equation}\label{eq:tangent}
 T(s)=\frac{A\cos(\omega s)+B\sin(\omega s)}
 {\abs{A\cos(\omega s)+B\sin(\omega s)}}.
\end{equation}
In this case the ground state is, up to a scalar factor,
\begin{equation}\label{eq:ground}
 f(s)=\abs{A\cos(\omega s)+B\sin(\omega s)}.
\end{equation}
\end{theorem}

The curvature operator already appeared in Burchard and Thomas's work on the dynamics of an inextensible elastic loop \cite{BT03}. Benguria and Loss connected its ground-state energy to the one-dimensional Lieb-Thirring inequality for two negative eigenvalues. For the original Lieb-Thirring conjecture and a broader account of the sharp-constant problem, see \cite{LT,FLW}. In \cite[Section~3]{BL}, a change to the cumulative mass variable converts the orthogonality and normalization of two eigenfunctions into the closure condition for a planar curve. Read and Schulz \cite{RS} have recently proved the sharp one-dimensional Lieb-Thirring inequality at exponent one, with constant $4/(3\sqrt3\pi)$, without a restriction on the number of bound states. The periodic geometric problem is distinct: the amplitude in the full-line transformation has vanishing endpoint values, while in \eqref{eq:rq} one must allow a positive periodic amplitude (see \cite[Remark~(i), Section~3]{BL}).

Benguria and Loss obtained the bound $1/2$ at length $2\pi$. Burchard and Thomas \cite{BT} established local minimality of the ovals. For convex planar curves, Linde \cite{Linde06} improved the lower bound to $(1+\pi/(\pi+8))^{-2}=0.6084\ldots$ and proved the conjecture under an additional geometric hypothesis. His later work \cite{Linde25} gives the analytic bound $0.81$ and a variational bound evaluated numerically as $0.8246$. The estimate
\[
 \lambda_1(\gamma)\geq \frac{\sqrt\pi}{2}
 \left(\frac{\Gamma(7/6)}{\Gamma(5/3)}\right)^3
 =0.9618\ldots
\]
in arbitrary dimension was obtained in \cite{Suragan} by a projection-averaging argument. The problem was discussed at the 2006 AIM workshop and the 2009 Oberwolfach workshop \cite{AIM,MFO}; see also the survey \cite[Section~8]{BLL}. Georgiev, G\'omez-Serrano, Tao and Wagner \cite[Section~8]{GGTW} used the conjecture as a benchmark for AlphaEvolve. Their numerical search recovered the circle but neither established the lower bound one nor recovered the other ovals.

Related questions concern the operator $-\partial_s^2+g\kappa^2$. Exner, Harrell and Loss \cite{EHL} proved circle minimality for planar curves when $0<g<1/4$ and showed that the circle need not minimize when $g>1$. The coefficient $g=1$ is special because of the family \eqref{eq:tangent}. Bernstein and Mettler \cite{BM} characterized this family by two geometric functionals and described its projective symmetries. The ovals conjecture also enters the study of minimal annuli in a slab \cite[Appendix~B]{BB}.

Denzler's existence and regularity theorem \cite[Theorem~1.1]{Denzler} reduces the global minimization problem, in arbitrary ambient dimension, to a strictly convex planar analytic curve with positive curvature. His analysis also excludes singular minimizers of the relaxed problem, apart from the digon \cite[Theorems~3.1 and~3.4]{Denzler}. The remaining difficulty is to understand the corresponding Euler equation. We address it by comparing the two arcs between antipodal tangent directions with harmonic oscillators having the same initial data. The comparison is not between their travel times: these may be different. Instead, it concerns an integral up to the first zero of a scalar solution (Lemma~\ref{lem:comparison}). 

In fact, comparison starting at both antipodal directions gives a rigidity statement for planar stationary curves of nonvanishing curvature, without an a priori bound on their eigenvalue (Theorem~\ref{thm:bias}). We conclude with the equality cases in Denzler's relaxed problem and a compactness consequence for almost minimizing sequences.

\section{A technical lemma}\label{sec:comparison}

For $a\in C([0,d])$, $a>0$, consider the system
\begin{equation}\label{eq:scalar}
 \dot u=ap,\quad \dot p=\left(\frac a4-\frac{3}{4a}\right)u.
\end{equation}
As usual, the dot denotes differentiation with respect to the time variable. Let $u(0)=u_*>0$, $p(0)=p_*$, and let $T_a$ be the first zero of $u$ in $(0,d]$, or $d$ if there is no such zero. We use the stopped integral
\[
 \mathcal M_d(a)=\int_0^{T_a} a(t)u(t)^2\dd t.
\]

\begin{lemma}\label{lem:comparison}
Assume $p_*/u_*\leq1/2$. If $a_1\geq a_0>0$ on $[0,d]$, then
\begin{equation}\label{eq:mono}
 \mathcal M_d(a_1)\geq\mathcal M_d(a_0).
\end{equation}
The inequality is strict if $a_1>a_0$ at any point of $[0,T_{a_0})$.
\end{lemma}

\begin{proof}
Put $b=a_1-a_0$ and $a_\eta=a_0+\eta b$, $0\leq\eta\leq1$. At a fixed value of $\eta$ we suppress its subscript on the solution and write $a=a_\eta$, $T=T_a$. With $X=(u,p)^\top$, the system has matrix
\[
 F(a)=\begin{pmatrix}0&a\\ a/4-3/(4a)&0\end{pmatrix},
 \quad
 F'(a)=\begin{pmatrix}0&1\\1/4+3/(4a^2)&0\end{pmatrix}.
\]
Solutions are continuously differentiable in $t$ and $\eta$, as follows from the integral equation for $X$ and the uniform positive lower bound for $a_\eta$.

We first justify differentiating the stopped integral. Every zero of $u$ is simple, since otherwise both components of $X$ vanish and uniqueness gives $X\equiv0$. If the first zero lies in $(0,d)$, the implicit function theorem continues it locally as a differentiable function of $\eta$. It remains the first zero by positivity on compact subintervals before it and the nonzero derivative at the zero. The boundary term in the derivative of $\mathcal M_d$ vanishes. If there is no zero in $[0,d]$, positivity persists under a small parameter variation. Finally, if the first zero is at $d$, extend the coefficients constantly just beyond $d$ and let $z(\eta)$ denote the continued zero. At the parameter $\eta_0$ where $z(\eta_0)=d$, one has $\abs{z(\eta)-d}=O(\abs{\eta-\eta_0})$. When $z(\eta)<d$, the difference between the integrals to $d$ and to $z(\eta)$ is $O(\abs{z(\eta)-d}^3)$, since $u_\eta$ vanishes at $z(\eta)$ and its time derivative is locally bounded. Thus the derivative is still obtained by differentiating the integrand. These arguments also give continuity of the stopped integral and of the stopped time.

Let $\Phi(s,t)$ be the fundamental matrix of $\dot X=F(a)X$, normalized by $\Phi(t,t)=I$, and put $e=(1,0)^\top$. Differentiation of the system and variation of constants give
\[
 \partial_\eta X(s)=\int_0^s\Phi(s,t)F'(a(t))X(t)b(t)\dd t.
\]
Substitution into the derivative of the stopped integral, followed by Fubini's theorem, yields
\begin{equation}\label{eq:kernel}
 \frac{\mathrm d}{\mathrm d\eta}\mathcal M_d(a_\eta)
 =\int_0^T b(t)K(t)\dd t,
\end{equation}
where
\begin{equation}\label{eq:K}
 K(t)=u(t)^2+2\int_t^T a(s)u(s)e^\top\Phi(s,t)F'(a(t))X(t)\dd s.
\end{equation}
We will show that $K(t)>0$ for $0\leq t<T$.

Fix such a $t$ and set $U=u(t)>0$, $L=p(t)/u(t)$ and
\[
 \sigma(s)=\int_t^s a(x)\dd x,\quad D=\sigma(T),\quad
 q(\sigma(s))=\frac{3}{4a(s)^2}.
\]
In this new time variable the first component of any homogeneous solution satisfies $y''=(1/4-q)y$, where the prime denotes differentiation in $\sigma$. Let $y,v$ be the solutions of this equation with initial data $(1,L)$ and $(0,1)$, respectively. Thus $u(s)=Uy(\sigma(s))$; in particular $y>0$ on $[0,D)$. The Wronskian identity gives
\begin{equation}\label{eq:wronskian}
 \frac{v(\sigma)}{y(\sigma)}=\int_0^\sigma\frac{\dd x}{y(x)^2},\quad 0\leq\sigma<D.
\end{equation}
It follows that $v>0$ on $(0,D)$. Moreover, on the positive part of the solution the slope $p/u$, in the time variable $\mathrm d\tau=a\,\mathrm dt$, satisfies
\[
 \frac{\mathrm d}{\mathrm d\tau}\left(\frac pu\right)
 =\frac14-\frac{3}{4a^2}-\left(\frac pu\right)^2.
\]
The right-hand side at $p/u=1/2$ is strictly negative. The initial slope condition consequently implies $L\leq1/2$.

At time $t$ the vector $F'(a(t))X(t)$ has components $UL$ and $U(1/4+3/(4a(t)^2))$. Propagating these data in the basis $y,v$ in the new time variable gives
\[
 e^\top\Phi(s,t)F'(a(t))X(t)
 =U\left[Ly(\sigma(s))+\left(\frac14-L^2+\frac{3}{4a(t)^2}\right)v(\sigma(s))\right].
\]
Using this identity and $\mathrm d\sigma=a(s)\,\mathrm ds$ in \eqref{eq:K}, we find
\begin{align}
 K(t)&=U^2\left(S+\frac{3}{2a(t)^2}\int_0^D yv\dd\sigma\right),\label{eq:kernel2}\\
 S&=1+2\int_0^D\left(Ly^2+\left(\frac14-L^2\right)yv\right)\dd\sigma.\notag
\end{align}
In particular the coefficient $a(t)$ in the second term is evaluated at the time of the variation, not at the integration variable.

To see that $S\geq0$, compare $y,v$ with the solutions of $Y''=Y/4$ having the same initial data:
\[
 Y_0(\sigma)=\cosh(\sigma/2)+2L\sinh(\sigma/2),\quad
 V_0(\sigma)=2\sinh(\sigma/2).
\]
Write $B_0=1/4-L^2$. We have $Y_0'=LY_0+B_0V_0$. Variation of constants gives
\[
 Y_0(\sigma)-y(\sigma)=2\int_0^\sigma
 \sinh((\sigma-x)/2)q(x)y(x)\dd x\geq0
\]
for $\sigma<D$, and the corresponding formula for $V_0-v$ gives $v\leq V_0$. Thus $0<y\leq Y_0$ and $0\leq v\leq V_0$ before $D$. The non-strict inequalities extend to $D$.

If $0\leq L\leq1/2$, the integrand defining $S$ is nonnegative. If $L\leq-1/2$, both $L$ and $B_0$ are nonpositive, and hence
\[
 Ly^2+B_0yv\geq LY_0^2+B_0Y_0V_0=Y_0Y_0'.
\]
Integration gives $S\geq Y_0(D)^2\geq0$.

It remains to consider $-1/2<L<0$. Now $B_0>0$ and $Y_0>0$ on $[0,\infty)$. By \eqref{eq:wronskian} and its counterpart for $Y_0,V_0$,
\[
 L+B_0\frac vy\geq L+B_0\frac{V_0}{Y_0}=\frac{Y_0'}{Y_0}.
\]
Where the left-hand side is negative, multiplication by $y^2\leq Y_0^2$ shows that $Ly^2+B_0yv\geq Y_0Y_0'$. Where it is nonnegative, the integrand is nonnegative. Consequently
\[
 S\geq1+2\int_0^\infty\min\{0,Y_0Y_0'\}\dd\sigma
 =\min_{\sigma\geq0} Y_0(\sigma)^2=1-4L^2>0.
\]
For the last equality, $Y_0'$ increases from $L<0$ to positive values and $(Y_0')^2-Y_0^2/4=L^2-1/4$.

In every case $S\geq0$. Since $D>0$ and $yv>0$ on $(0,D)$, \eqref{eq:kernel2} gives strict positivity of $K$ on $[0,T)$. The expression \eqref{eq:K} is continuous up to the stopped endpoint. Equation \eqref{eq:kernel} now proves monotonicity. If $b>0$ at a point before $T_{a_0}$, it is positive on a subinterval there. By continuity of the stopped time, that subinterval (possibly shortened) is still before $T_{a_\eta}$ for $\eta>0$ sufficiently small. The derivative is strictly positive for such $\eta$, which proves strictness in \eqref{eq:mono}.
\end{proof}

\section{Euler equation and its two half-turns}\label{sec:ode}

By scaling it suffices to prove Theorem~\ref{thm:main} at length $2\pi$. For a positive amplitude $f$ put $\Psi=fT$. The basic observation of \cite{BL,BT} (see also \cite[Lemma~2.2]{Denzler}) is
\begin{equation}\label{eq:vector}
 \int_{\T_{2\pi}}\abs{\Psi'}^2\dd s
 =\int_{\T_{2\pi}}(f'^2+\kappa^2f^2)\dd s,
 \quad \mathcal C(\Psi):=\int_{\T_{2\pi}}\frac{\Psi}{\abs{\Psi}}\dd s=0.
\end{equation}
Conversely, every nowhere-zero $H^1$ vector field satisfying the latter condition defines a closed $W^{2,2}$ curve by integrating its unit direction.

Denzler \cite[Theorem~1.1]{Denzler} proved that the infimum of \eqref{eq:rq} is attained and that every minimizing curve is planar, strictly convex and analytic, with positive curvature. Its tangent makes exactly one turn. Let $f>0$ be the normalized ground state of a minimizing curve, and write $\lambda=\lambda_1(\gamma)$. The ground state can be chosen positive: a nonnegative solution of $f''=(\kappa^2-\lambda)f$ cannot vanish without its derivative vanishing as well, contrary to uniqueness. The circle gives $\lambda\leq1$, and $\lambda>0$ since zero energy would require a nonzero constant amplitude and zero curvature, contradicting closure.

For the constrained vector quotient \eqref{eq:vector}, the derivative of $\mathcal C$ is
\[
 D\mathcal C(\Psi)V=\int_{\T_{2\pi}}\frac{(I-TT^\top)V}{f}\dd s.
\]
It is onto: a vector annihilating its range would have to be parallel to every tangent, which is impossible for a closed regular curve. The mass derivative is independent of this constraint, as is seen by taking $V=\Psi$. Lagrange multipliers therefore yield
\begin{equation}\label{eq:EL}
 \Psi''+\lambda\Psi=-\frac{1}{\abs{\Psi}}(I-TT^\top)\mu,
 \quad \mu\in\R^n.
\end{equation}
Here the coefficient of $\Psi$ is identified by integration against $\Psi$. Since $\Psi$ lies in a plane, projection onto its orthogonal complement shows that $\mu$ lies in that plane as well. Equation \eqref{eq:EL} is the stationary equation studied in \cite{BT} and \cite[Equation~(3.2)]{Denzler}.

Set $t=\sqrt\lambda\,s$, $Y(t)=\Psi(s)$ and $P=2\pi\sqrt\lambda$. If $\mu\ne0$, orient the first coordinate axis along $\mu$ and put $c=\abs\mu/\lambda$. Then
\begin{equation}\label{eq:force}
 \ddot Y+Y=-c\nabla H(Y),\quad H(Y)=\frac{Y_1}{\abs{Y}}.
\end{equation}
The energy is conserved (cf. \cite[Equation~(3.3)]{Denzler}):
\begin{equation}\label{eq:energy}
 \abs{Y}^2+\abs{\dot Y}^2=E-2cH(Y).
\end{equation}
The closure condition is $\int_0^P Y/\abs Y\,\mathrm dt=0$. We will analyze \eqref{eq:force} for either sign of $c$ and without imposing $P\leq2\pi$.

\subsection*{Radial comparison}

Consider any nonvanishing planar solution of \eqref{eq:force} whose polar angle is strictly monotone with nonzero angular velocity. Start at the direction $-e_1$ and travel in either direction in time until the first arrival at $e_1$, say at time $h$. Time reversal leaves \eqref{eq:force} unchanged. Each of these half-turns has a continuous angle $\theta$ from $\pi$ to $0$ or from $\pi$ to $2\pi$. In particular, $H(Y)=\cos\theta$ increases strictly from $-1$ to $1$ on each outward arc.

At the initial point, $E+2c=\abs{Y(0)}^2+\abs{\dot Y(0)}^2>0$. Scaling $Y$ by a positive constant rescales $E,c$ by its square, since $H$ is homogeneous of degree zero. We may therefore assume $E+2c=1$, with the sign of $c$ unchanged. Put $r=\abs Y$ and compare $Y$ with
\[
 Q(t)=Y(0)\cos t+\dot Y(0)\sin t,\quad R(t)=\abs{Q(t)}.
\]
The two initial vectors are linearly independent because the angular velocity is nonzero. Thus $R>0$ everywhere, and the first arrival of $Q/\abs Q$ at $e_1$ is at $t=\pi$.

By homogeneity $Y\cdot\nabla H(Y)=0$. Equations \eqref{eq:force}--\eqref{eq:energy} imply
\[
 \frac{\mathrm d^2}{\mathrm dt^2}r^2+4r^2=2(E-2cH(Y)),
 \quad \frac{\mathrm d^2}{\mathrm dt^2}R^2+4R^2=2.
\]
The two squared radii have identical initial values and first derivatives. Solving the difference equation and integrating by parts, using $1+H(Y(0))=0$, gives
\begin{align}
 R(t)^2-r(t)^2
 &=2c\int_0^t\sin(2(t-x))(1+H(Y(x)))\dd x\notag\\
 &=2c\int_0^t\sin^2(t-x)\frac{\mathrm d}{\mathrm dx}H(Y(x))\dd x.
 \label{eq:radial}
\end{align}
For $0<t\leq h$ the integral on the last line is positive. Hence $r<R$ when $c>0$, and $r>R$ when $c<0$. The comparison holds up to the actual arrival time, which need not be $\pi$.

\subsection*{The square-root coordinate}

On an outward arc write $Y=r(\cos\theta,\sin\theta)$ and set
\[
 \mathrm d\tau=\frac{\mathrm dt}{r},\quad u=\sqrt r\sin(\theta/2),\quad p=\frac{\mathrm du}{\mathrm d\tau}.
\]
This is one component of the classical square-root coordinate associated with Levi-Civita's regularization \cite{LC}. With our choice of angle, $u>0$ until the arrival at $e_1$, where it has its first zero. Set $v=\dot r$, $z=r\dot\theta$ and $w=\theta/2$. The radial and angular equations from \eqref{eq:force} give
\[
 r_\tau=rv,\quad v_\tau=z^2-r^2,\quad
 \theta_\tau=z,\quad z_\tau=c\sin\theta-vz.
\]
In particular $2p=\sqrt r(v\sin w+z\cos w)$. Differentiating this expression and using the energy identity yields
\begin{align*}
 4u_{\tau\tau}
 &=\sqrt r\big((v^2+z^2-2r^2)\sin w+2c\sin\theta\cos w\big)\\
 &=(E+2c-3r^2)u=(1-3r^2)u.
\end{align*}
Here we used $2\sin\theta\cos w=2(1+\cos\theta)\sin w$. Returning to physical time gives precisely \eqref{eq:scalar} with $a=1/r$. The initial data are
\begin{equation}\label{eq:initial}
 u_* =\sqrt{r(0)},\quad \frac{p_*}{u_*}=\frac{\dot r(0)}2<\frac12,
\end{equation}
since $r(0)^2+\abs{\dot Y(0)}^2=1$. The free reference $Q$ satisfies the same scalar equation with $a_0=1/R$ and the same initial data. Its corresponding first zero is at $\pi$. Notice also that the stopped mass has the simple interpretation
\begin{equation}\label{eq:mass}
 M:=\int_0^h au^2\dd t=\frac12\int_0^h (1-H(Y(t)))\dd t.
\end{equation}
Let $M_0$ be the corresponding integral for $Q$ on $[0,\pi]$.

\begin{proposition}\label{prop:arc}
For each half-turn starting at $-e_1$, $M>M_0$ if $c>0$, and $M<M_0$ if $c<0$.
\end{proposition}
\begin{proof}
Use a common horizon $d=\max\{h,\pi\}$. If $h<\pi$, continue the coefficient $a$ by
\[
 a(t)=a(h)\frac{a_0(t)}{a_0(h)},\quad h\leq t\leq\pi.
\]
This does not change the scalar solution before its first zero at $h$, and it preserves the coefficient ordering. If $h\geq\pi$, no extension is necessary; the reference $Q$ is already defined globally. In either case $M=\mathcal M_d(a)$ and $M_0=\mathcal M_d(a_0)$. The radial identity \eqref{eq:radial} implies $a>a_0$ on $(0,h]$ for $c>0$, and the reverse inequality for $c<0$. Both coefficients have the data \eqref{eq:initial}. Lemma~\ref{lem:comparison}, with its strictness assertion, proves the claim for either ordering.
\end{proof}

\section{Rigidity from the closure condition}\label{sec:rigidity}

We first record what the half-turn comparison says about periodic orbits of the forced oscillator itself. 

\begin{theorem}\label{thm:bias}
Let $Y\colon\R/P\mathbb Z\to\R^2\setminus\{0\}$ solve \eqref{eq:force}, with $c\ne0$. Suppose that its polar angle has nonzero angular velocity of constant sign and makes $m$ full turns in absolute value in the period $P$, where $m$ is a positive integer. Then
\begin{equation}\label{eq:bias}
 \operatorname{sgn}(c)\int_0^P H(Y(t))\dd t<-\abs{P-2m\pi}.
\end{equation}
In particular, its normalized direction cannot have zero mean.
\end{theorem}
\begin{proof}
Suppose first that $c>0$ and write $J=\int_0^P H(Y(t))\dd t$. At each of the $m$ occurrences of the direction $-e_1$, form the two outward half-turns used in Section~\ref{sec:ode}. For each such pair the initial velocities are negatives of each other. If the harmonic reference on one arc is $Q(t)$, it is $Q(-t)$ on the other. Hence the sum of the two reference masses for each pair is
\[
 M_{0,+}+M_{0,-}=\frac12\int_{-\pi}^{\pi}(1-H(Q(t)))\dd t=\pi,
\]
since $Q(t+\pi)=-Q(t)$. The $2m$ actual arcs cover the period exactly once, with reversed time on half of them, and their masses sum to $(P-J)/2$. The reference masses sum to $m\pi$. Proposition~\ref{prop:arc} yields
\begin{equation}\label{eq:firstbias}
 P-J>2m\pi.
\end{equation}

Now apply the same argument to $-Y$, which satisfies \eqref{eq:force} with coefficient $-c$ because $\nabla H(-Y)=\nabla H(Y)$. The normalization of energy may be different, but it does not affect $H$, the period, or the masses in \eqref{eq:mass}. Since $H(-Y)=-H(Y)$, the reverse comparison in Proposition~\ref{prop:arc} gives
\begin{equation}\label{eq:secondbias}
 P+J<2m\pi.
\end{equation}
Together \eqref{eq:firstbias} and \eqref{eq:secondbias} imply $J<\min\{P-2m\pi,2m\pi-P\}=-\abs{P-2m\pi}$. The case $c<0$ follows by applying this result to $-Y$.
\end{proof}

\begin{proof}[Proof of Theorem~\ref{thm:main}]
Return to the minimizing curve at length $2\pi$ supplied by Denzler's theorem. If its multiplier $\mu$ in \eqref{eq:EL} were nonzero, the rescaled orbit would satisfy Theorem~\ref{thm:bias}, contradicting closure. Thus $\mu=0$ and
\[
 \Psi(s)=A\cos(\sqrt\lambda\,s)+B\sin(\sqrt\lambda\,s).
\]
The vectors $A,B$ are linearly independent because $\Psi$ never vanishes. Periodicity implies $\sqrt\lambda$ is a positive integer, and $\lambda\leq1$ forces $\lambda=1$. This establishes the lower bound. Every equality curve is itself a global minimizer, so the same argument gives the stated form of its tangent and ground state.

For the converse, at any length $\ell$, take independent $A,B$, put $Y=A\cos(\omega s)+B\sin(\omega s)$, and set $f=\abs Y$, $T=Y/f$. The relation $T(s+\ell/2)=-T(s)$ guarantees closure. The identity $Y''+\omega^2Y=0$, projected onto $T$, gives $-f''+\kappa^2 f=\omega^2 f$. Since $f>0$, the ground-state identity
\begin{equation}\label{eq:groundstate}
 \int_{\T_\ell}\big(\abs{g'}^2+(\kappa^2-\omega^2)\abs g^2\big)\dd s
 =\int_{\T_\ell}f^2\left|\left(\frac gf\right)'\right|^2\dd s
\end{equation}
holds for every real or complex $g\in H^1(\T_\ell)$, by integration by parts. Thus the first eigenvalue is $\omega^2$, and the equality functions are exactly the scalar multiples of $f$. Finally, dilation by $\omega$ changes length $\ell$ to $2\pi$ and multiplies eigenvalues by $\omega^{-2}$, giving \eqref{eq:main} at general length.
\end{proof}

At length $2\pi$, an arclength shift diagonalizes the Gram matrix of $A,B$. After a rotation or reflection and multiplication of the pair by a common positive scalar, the tangent of an equality curve is
\[
 T_\alpha(s)=\frac{(\alpha\cos s,\alpha^{-1}\sin s)}
 {\sqrt{\alpha^2\cos^2s+\alpha^{-2}\sin^2s}},\quad \alpha\geq1.
\]
Its curvature is $(\alpha^2\cos^2s+\alpha^{-2}\sin^2s)^{-1}$, as follows by computing the determinant of $Y,Y'$. The curve is obtained by integration of $T_\alpha$. The value $\alpha=1$ gives the unit circle, while as $\alpha\to\infty$ the tangents converge almost everywhere to a direction which reverses twice per period. The limit is a segment traversed twice and does not belong to $W^{2,2}$.

\begin{corollary}\label{cor:stationary}
Let $\gamma$ be a smooth closed planar curve of length $\ell$ with nowhere-vanishing curvature whose tangent makes $m$ turns in absolute value. It is stationary for the lowest periodic curvature eigenvalue at fixed length if and only if it is an $m$-fold traversal of an oval of length $\ell/m$. In this case $\lambda_1(\gamma)=(2m\pi/\ell)^2$. In particular every strictly convex stationary curve of positive curvature is a global minimizer.
\end{corollary}
\begin{proof}
At length $2\pi$, take the positive ground state $f$ of such a curve. Its eigenvalue $\lambda$ is positive. The vector $\Psi=fT$ is stationary for the quotient in \eqref{eq:vector} subject to closure. Indeed, along any smooth admissible variation of $\Psi$, the quotient is at least the eigenvalue of the corresponding curve and agrees with it at the initial point. The eigenvalue has zero first variation by hypothesis, so the quotient has zero first variation as well. The argument leading to \eqref{eq:EL} therefore applies. Theorem~\ref{thm:bias} forces $\mu=0$, without needing $\lambda\leq1$. The angle of $\Psi=A\cos(\sqrt\lambda\,s)+B\sin(\sqrt\lambda\,s)$ winds $\sqrt\lambda$ times in the positive or negative orientation in length $2\pi$. The turning number in absolute value is $m$, so $\lambda=m^2$. The tangent is the normalization of this harmonic orbit. It is periodic with period $2\pi/m$ and changes sign after half that period. Thus its integral on each period vanishes, and the curve traverses the corresponding oval $m$ times. Rescaling gives the assertion at length $\ell$. Conversely, on an $m$-fold oval the positive amplitude in \eqref{eq:ground}, with frequency $2m\pi/\ell$ instead of $\omega$, is a ground state by the same argument as in \eqref{eq:groundstate}. The corresponding vector field is harmonic at that frequency, so its vector Rayleigh quotient has zero first variation, even without imposing closure. Along any smooth fixed-length variation of the curve, use this fixed amplitude as a test function. Its Rayleigh quotient is stationary and touches the lowest eigenvalue from above at the original curve. The lowest eigenvalue is differentiable along the variation (it is simple), so its first derivative also vanishes. Finally, $m=1$ gives global minimality by Theorem~\ref{thm:main}.
\end{proof}

\section{Relaxed problem and almost minimizers}\label{sec:relaxed}

In the vector formulation it is natural to allow zeros. For nonzero $\Psi\in H^1(\T_{2\pi};\R^n)$ let $Z=\{s:\Psi(s)=0\}$. Denzler's weak loop constraint \cite[Definition~2.3, Equation~(2.6)]{Denzler} is
\begin{equation}\label{eq:weakclosure}
 \left|\int_{\T_{2\pi}\setminus Z}\frac{\Psi}{\abs\Psi}\dd s\right|\leq |Z|,
\end{equation}
where $|Z|$ denotes Lebesgue measure. This permits the tangent direction to be filled in on the zero set. Define $\mathcal E$ to be the family of vector fields of $L^2$ norm one described in the next theorem.

\begin{theorem}\label{thm:relaxed}
If $\Psi$ satisfies \eqref{eq:weakclosure}, then
\begin{equation}\label{eq:weakineq}
 \int_{\T_{2\pi}}\abs{\Psi'}^2\dd s\geq\int_{\T_{2\pi}}\abs\Psi^2\dd s.
\end{equation}
Equality holds precisely for the following fields, up to a translation of the parameter:
\begin{enumerate}
 \item $\Psi=A\cos s+B\sin s$, with $A,B\in\R^n$ linearly independent;
 \item $\Psi$ of the form
 \begin{equation}\label{eq:degenerate}
 \Psi(s)=\begin{cases}
 a e\sin s,&0\leq s\leq\pi,\\
 b e\sin s,&\pi\leq s\leq2\pi,
 \end{cases}
 \quad \abs e=1,\quad a,b\geq0,\quad a+b>0.
 \end{equation}
\end{enumerate}
\end{theorem}
\begin{proof}
By \cite[Theorems~3.1 and~3.4]{Denzler}, the relaxed quotient attains its minimum, and every minimizer has either no zeros, one active interval of length $\pi$, or two active intervals each of length $\pi$. In the first case its modulus and direction define an ordinary curve, and Theorem~\ref{thm:main} applies. In the other two cases the Dirichlet Poincar\'e inequality on each active interval gives a quotient at least one. The circle attains one, so \eqref{eq:weakineq} follows.

An equality field is a relaxed minimizer and thus satisfies the same alternative for the zero set. If the field is nowhere zero, its modulus must be a ground state of its associated curve. Theorem~\ref{thm:main} gives the first family. Otherwise equality in Dirichlet Poincar\'e forces a constant vector times the first sine mode on each active interval. When both halves are active, their unit directions must be opposite by \eqref{eq:weakclosure}. When only one is active, the zero set has measure $\pi$ and the constraint allows an arbitrary fixed direction on the active half. This gives \eqref{eq:degenerate}. Direct substitution verifies that both families are admissible and attain equality.
\end{proof}

There is a compactness consequence which takes the degenerating ovals into account. It concerns the amplitude-weighted tangent, not the tangent alone.

\begin{corollary}\label{cor:compactness}
For each fixed $n\geq2$ and every $\varepsilon>0$ there is $\delta>0$ such that every $\Psi\in H^1(\T_{2\pi};\R^n)$ satisfying \eqref{eq:weakclosure}, $\norm\Psi_{L^2}=1$, and
\[
 \int_{\T_{2\pi}}\abs{\Psi'}^2\dd s\leq1+\delta
\]
satisfies $\operatorname{dist}_{H^1}(\Psi,\mathcal E)<\varepsilon$.
\end{corollary}
\begin{proof}
If the assertion fails, take a sequence $\Psi_j$ satisfying the normalization and weak loop constraint, with energies tending to one but with distance from $\mathcal E$ bounded below by a positive constant. The $H^1$ bound gives a subsequence converging uniformly and weakly in $H^1$ to $\Psi$, with $\norm\Psi_{L^2}=1$. To check the constraint, let $Z_j=\{\Psi_j=0\}$ and $Z=\{\Psi=0\}$, defining $\operatorname{sgn}(0)=0$ for vector fields. Then
\[
 \left|\int_{\T_{2\pi}\setminus Z}\operatorname{sgn}\Psi_j\dd s\right|
 \leq |Z_j|+|Z\setminus Z_j|=|Z|+|Z_j\setminus Z|.
\]
On the complement of $Z$, the signs converge pointwise and $\mathbf1_{Z_j}$ tends to zero. Dominated convergence gives \eqref{eq:weakclosure} for $\Psi$. By lower semicontinuity and Theorem~\ref{thm:relaxed}, its energy is one. Consequently the derivative norms converge, and weak $H^1$ convergence improves to strong $H^1$ convergence. The limit belongs to $\mathcal E$, contradicting the distance bound.
\end{proof}

\newpage
{\bf Funding.} This research is funded by Nazarbayev University under the FDCRGP 110326FD3205 (D.S.).

\medskip
{\bf Use of AI.} The author acknowledges the use of several AI tools during the study. The main
mathematical arguments are obtained by the author. OpenAI ChatGPT assisted in revising/writing the final version.

\end{document}